\magnification=\magstep1

\input epsf

~~\vskip .1in

~~\vskip 1in

\centerline{\bf Probing mapping class groups using
arcs}

\vskip .2in

\centerline {\bf R. C. Penner}

\centerline{Departments of Mathematics and Physics/Astronomy}

\centerline{University of Southern California}

\centerline {Los Angeles, CA 90089}

\vskip .2in

\leftskip .8in\rightskip .8in

\noindent {\bf Abstract}~~The action of the mapping class
group of a surface on the collection of homotopy classes of
disjointly embedded curves or arcs in the surface is discussed here as a
tool for understanding Riemann's moduli space and its topological and geometric
invariants.  Furthermore, appropriate completions, elaborations, or quotients of the set of
all such homotopy classes of curves or arcs give for instance Thurston's boundary for
Teichm\"uller space or a combinatorial description of moduli
space in terms of fatgraphs.  Related open problems and questions are discussed.

\vskip .3in

\leftskip=0ex\rightskip=0ex

\centerline{\bf Introduction}

\vskip .2in

\noindent One basic
theme in this paper on open problems is that the action of
the mapping class group on spaces of measured foliations, and in particular on weighted families
of curves and arcs, is calculable and captures the
dynamics of homeomorphisms of the surface
both on the surface itself and on its Teichm\"uller space.  Another basic theme is that 
suitable spaces of arcs can be exploited to give group-theoretic and other data about the mapping class 
groups as well as their subgroups and completions.
The author was specifically given the task by the editor of presenting
open problems on his
earliest and latest works in [23] and [46-49], which respectively develop these two themes and are
also surveyed here.

\vskip .3in

\noindent {\bf 1. Definitions}

\vskip .2in
    
\noindent Let $F=F_{g,r}^s$ denote a smooth oriented surface of genus
$g\geq 0$ with $r\geq 0$ boundary components and $s\geq 0$ punctures, where
$2g-2+r+s> 0$.  The mapping class group $MC(F)$ of $F$ is the collection of isotopy
classes of all orientation-preserving homeomorphisms of $F$, where the isotopies and
homeomorphisms necessarily setwise fix the boundary $\partial F$ and setwise fix the collection of punctures.  Let $PMC(F)<MC(F)$
denote the pure mapping class group whose homeomorphisms and isotopies pointwise fix each puncture and each boundary component.

\vskip .1in

\noindent But one aspect of Bill Thurston's seminal contributions to mathematics, [37] among other things
provides a natural spherical boundary of the Teichm\"uller space ${\cal T}(F)$ of a surface $F$ with negative Euler characteristic
by an
appropriate space of ``projectivized measured foliations of compact support'' in $F$. 
Specifically [37,9,32], let ${\cal MF}(F)$ denote the space of all
isotopy classes rel $\partial F$ of Whitehead equivalence classes of measured foliations in
$F$, and let ${\cal
MF}_0(F)$ denote the subspace comprised of those foliations with ``compact support'',
i.e., leaves are disjoint from a neighborhood of the punctures and boundary, and no simple closed leaf is puncture-
or boundary-parallel.

\vskip .1in

\noindent Define the projectivized spaces
$$\eqalign{
{\cal PF}_0(F)~&=~[{\cal MF}_0(F)-\{ \vec 0\} ]/{\bf R_{>0}}\cr
&\subseteq
~~[{\cal MF}(F)-\{ \vec 0\} ]/{\bf R_{>0}}~=~{\cal PF}(F),\cr
}$$
where $\vec 0$ denotes the empty foliation and ${\bf R}_{>0}$ acts by homothety on
transverse measures.  Thus, 
${\cal MF}(F_{g,r-1}^{s+1})\subseteq{\cal
MF}(F_{g,r}^{s})$ and
${\cal MF}_0(F_{g,r}^s)\approx {\cal
MF}_0(F_{g,0}^{r+s})$ with corresponding statements also for projectivized foliations.

\vskip .1in

\noindent
A basic fact (as follows from the density of simple closed curves in ${\cal PF}(F)$) is that the action of $MC(F)$ or $PMC(F)$ on ${\cal
PF}_0(F)$ has dense orbits, so the quotients are
non-Hausdorff.  (The action is actually minimal in the sense that every orbit is dense, and in fact, the action is ergodic for a 
natural measure class as independently shown by
Veech [43] and Masur[44].)  

\vskip .1in

\noindent We shall say that a measured foliation or its projective class {\it fills} $F$ if every essential simple
closed curve has positive transverse
measure and that it {\it quasi fills} $F$ if every essential curve with
vanishing transverse measure is puncture-parallel.

\vskip .1in

\noindent Define the {\it pre-arc complex} ${\cal A}'(F)$ to be
the subspace of
${\cal PF}(F)$ where each leaf in the underlying foliation is required to be an arc
connecting punctures or boundary components, and define the open subspace ${\cal
A}_\#'(F)\subseteq{\cal A}'(F)$ where the foliations are furthermore required to fill $F$.   

\vskip .1in

\noindent In
particular in the {\it punctured case} when $r=0$,
$s\geq 1$, and $F$ has negative Euler characteristic the product
${\cal T}(F)\times \Delta ^{s-1}$ of Teichm\"uller space
with an open 
($s-1$)-dimensional simplex $\Delta ^{s-1}$ is
canonically isomorphic to ${\cal A}'_\# (F)$, and in fact this descends to an isomorphism
between the {\it filling arc complex} ${\cal A}_\#(F)={\cal
A}_\#'(F)/PMC(F)$ and the product
${\cal M}(F)\times \Delta ^{s-1}$ of Riemann's (pure) moduli space ${\cal M}(F)={\cal
T}(F)/PMC(F)$ with the simplex [12,15,36,27] .  Thus, the {\it arc complex} ${\cal A}(F)={\cal A}'(F)/PMC(F)$ forms a natural
combinatorial compactification of ${\cal A}_\#(F)\approx{\cal M}(F)\times\Delta^{s-1}$.
(In fact, this is {\sl not} the most useful combinatorial compactification when $s>1$, cf. [30,46], where one chooses from among the punctures a
distinguished one.)

\vskip .1in

\noindent Another special case $r,s$ of interest here is the case of {\it bordered surfaces}
when $r\geq 1$ and $s\geq 0$.  Choose one distinguished point on each boundary component, and define the analogous complexes ${Arc}'(F)
\subseteq {\cal PF}(F)$, where leaves are required to be asymptotic to the distinguished points on the boundary
(and may not be asymptotic to punctures)
with its quasi filling subspace ${Arc}'_\#(F)\subseteq	{Arc}'(F)$ and
quotients
${Arc}_\#(F)={Arc}_\#'(F)/PMC(F)\subseteq {Arc}(F)={Arc}'(F)/PMC(F)$. 
In analogy to the punctured case, 
${Arc}(F)$ is proper homotopy equivalent to Riemann's moduli space of $F$ (as a bordered
surface with one distinguished point in each geodesic boundary component) modulo a natural action of the positive reals [31].
Furthermore [30], $Arc$-complexes occur as virtual links of simplices in the ${\cal A}$-complexes, and
the local structure of the ${\cal A}$-complexes is thus governed by the global topology of
${Arc}$-complexes.  In fact, the ${Arc}$-complexes are stratified spaces of a particular sort as explained
in
$\S$5.

\vskip .1in

\noindent There are other geometrically interesting subspaces and quotients of ${\cal
MF}(F)$ or ${\cal PF}(F)$, for instance the curve complex of Harvey [14] or the complex of pants
decompositions of Hatcher-Thurston [13], which are surely discussed elsewhere in this volume.

\vskip .1in

\noindent As was noted before, the quotients ${\cal PF}_0(F)/PMC(F)\subseteq {\cal
PF}(F)/PMC(F)$ are maximally non-Hausdorff, and yet for $r>0$, ${\cal PF}(F)/PMC(F)$ contains
as an open dense subset the (Hausdorff) stratified space ${Arc}(F)=Arc'(F)/PMC(F)$; in particular, for the
surface $F=F_{0,r}^s$ with $r+s\leq 3$, 
$PMC(F_{0,r}^s)$ is the free abelian group generated by Dehn twists
on the boundary, ${Arc}(F_{0,r}^s)$ is piecewise-linearly homeomorphic to a sphere of 
dimension
$3r+2s-7$, and it is not difficult to understand the non-Hausdorff space
${\cal PF}(F)/PMC(F)$ with its natural foliation.
This leads to
our first problem:
 
\vskip .2in

\noindent {\bf Problem 1}~~Understand either classically or as quantum geometric objects
the non-Hausdorff quotients of ${\cal PL}_0(F)$ or ${\cal PL}(F)$ by $MC(F)$ or $PMC(F)$.

\vskip .2in

\noindent ${\cal T}(F)$ has been quantized in [5] and  [17] as surveyed
in [6] and [35], respectively, and 
${\cal PL}_0(F_{1,0}^1)$ has been quantized in [6].  In interesting contrast,
[21,54] has described a program for studying real quadratic number fields as
quantum tori limits of elliptic curves, thus quantizing limiting curves rather than Teichm\"uller space.

\vskip .3in

\noindent{\bf 2. Dehn-Thurston coordinates}

\vskip .2in

\noindent Fix a surface $F=F_{g,r}^s$.  In this section, we introduce global
Dehn-Thurston coordinates from [23] for ${\cal MF}_0(F)$ and ${\cal MF}(F)$.

\vskip .1in

\noindent  Define a {\it pants decomposition} of $F$ to be a
collection ${\cal P}$ of curves disjointly embedded in $F$ so that each component of
$F-\cup{\cal P}$ is a {\it generalized pair of pants} $F_{0,r}^s$ with $r+s=3$.  One easily
checks using Euler characteristics that there are $3g-3+2r+s$ curves in a pants
decomposition of $F$ and $2g-2+r+s$ generalized pairs of pants complementary to ${\cal P}$.

\vskip .1in

\noindent Given a measured foliation ${\cal F}$ in a generalized pair of pants $P$, where $\partial P$ has components $\partial _i$,
define the triple $m_i$ of {\it intersection numbers} given by the transverse measure of $\partial _i$ of
${\cal F}$, for $i=1,2,3$.

\vskip .2in

\noindent {\bf Dehn-Thurston Lemma}~\it Isotopy classes (not necessarily the identity on the boundary) of Whitehead equivalence classes of
non-trivial measured foliations in the pair of pants $P$ with no closed leaves are uniquely determined by the triple $(m_1,m_2,m_3)$ of
non-negative real intersection numbers, which are subject only to the constraint that
$m_1+m_2+m_3$ is positive.  \rm

\vskip.2in

\vskip .2in

~~~~\epsffile{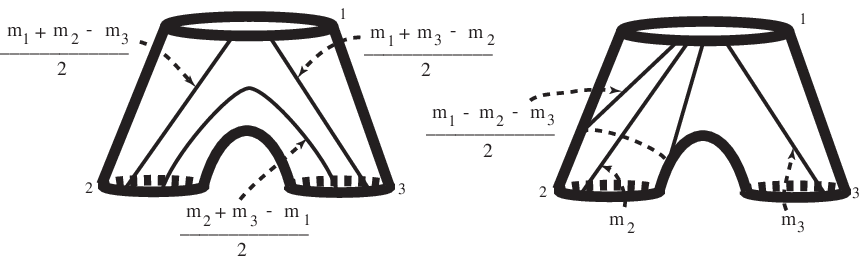}

\hskip .4in{\bf 1a}~Triangle inequality\hskip .5in{\bf 1b}~The case $m_1>m_2+m_3$

\vskip .1in

\centerline{{\bf Figure 1}~Constructing measured foliations in pants}

\vskip .1in

\noindent {\bf Proof}~~The explicit construction of a measured foliation realizing a putative triple of intersection numbers
is illustrated in Figure~1 in the two representative cases that $m_1,m_2,m_3$ satisfy all three possible (weak) triangle inequalities
in Figure~1a
or perhaps $m_1\geq m_2+m_3$ in Figure~1b.  The
other cases are similar, and the projectivization of the positive orthant in
$(m_1,m_2,m_3)$-space is illustrated in Figure~2.  Elementary topological considerations show that any measured foliation of $P$ is isotopic
to a unique such foliations keeping endpoints of arcs in the boundary of $P$, completing the proof.~~~~\hfill{\it q.e.d.}

\vskip .1in

\centerline{\epsffile{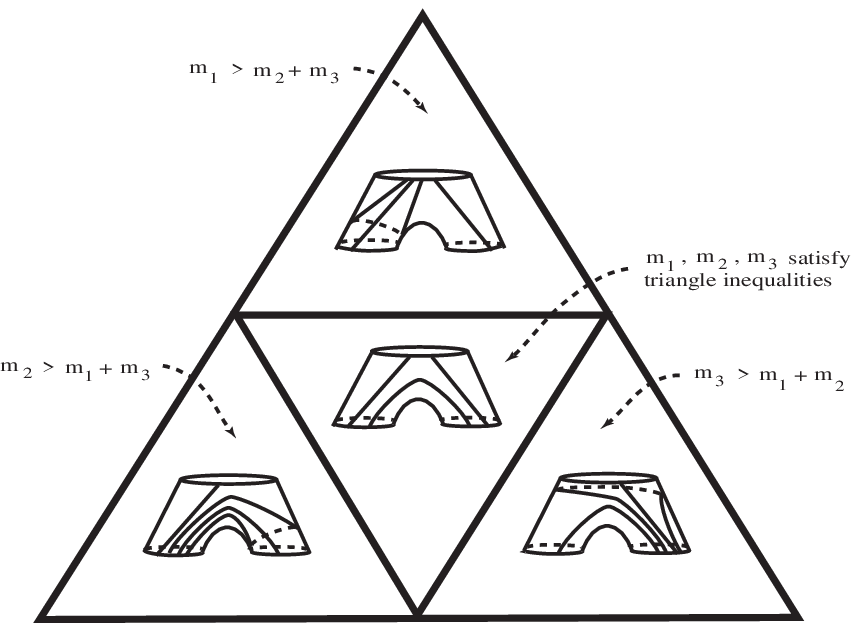}}

\vskip .1in

\centerline{{\bf Figure 2}~Measured foliations in pants}

\vskip .1in

\noindent In order to refine the Dehn-Thurston Lemma and keep track of twisting around the boundary, we shall introduce in each
component
$\partial _i$ of the boundary $\partial F$ an arc $w_i\subseteq \partial _i$ called a {\it window}, for
$i=1,2,3$.  We require that the support of the restriction to $\partial P$  of ${\cal F}$ lie in the union of
the windows, so-called {\it windowed measured foliations}.
(Collapsing each window to a point gives a surface with a distinguished point in each boundary component, so a windowed measured foliation in $P$
gives
rise  to an element of ${Arc}(P)$ in the sense of $\S$1.)  We seek the analogue of the Dehn-Thurston Lemma for
windowed isotopy classes.

\vskip .1in

\noindent To this end, there are two conventions to be made:

\vskip .1in

\leftskip .3in

\noindent 1) when a leaf of ${\cal F}$
connects a boundary component to itself, i.e., when it is a loop, then it passes around a specified leg, right or left, of
$P$ as illustrated in Figure 3a-b, i.e., it contains a particular boundary component or puncture in its complementary component with one
cusp;

\vskip .1in 

\noindent 2) when a leaf is added to the complementary component of a loop in $P$ with two cusps, then it either
follows or precedes the loop as illustrated in Figure 3c-d.

\leftskip=0ex

\vskip .2in 

\noindent For instance in Figure~2, the conventions are: 1) around the right leg
for loops; and 2) on the front of the surface.  We shall call these the {\it standard twisting conventions}.

\vskip .1in

\hskip .1in{{{\epsffile{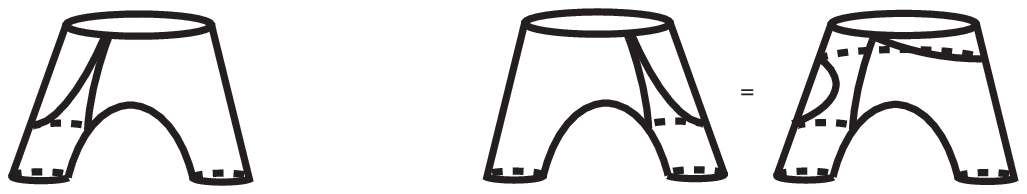}}}}

{\bf Figure 3a} Around the right leg.\hskip .4in{\bf Figure
3b} Around the left leg.

\vskip .05in

\hskip .1in{{{\epsffile{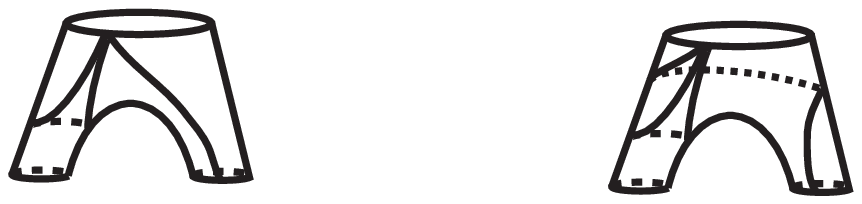}}}}

{\bf Figure 3c} Arc follows loop.\hskip .5in~~~~~~{\bf Figure
3d}~Arc precedes loop.

\vskip .1in

\centerline{{\bf Figure 3}~~{Twisting conventions.}}

\vskip .1in

\noindent Upon making such choices of convention, we may associate a twisting number $t_i\in {\bf R}$ to ${\cal F}$ as follows.  Choose a regular
neighborhood of $\partial P$ and consider the sub-pair of pants
$P_1\subseteq P$ complementary to this regular neighborhood.  Given a weighted arc family $\alpha$ in $P$, by the Dehn-Thurston Lemma,
we may perform an isotopy in $P$ supported on a neighborhood of $P_1$ to arrange that $\alpha\cap P_1$ agrees with a conventional
windowed arc family in $P_1$ (where the window in $\partial P_1$ arises from that in $\partial P$ in the natural way via a framing of the
normal bundle to $\partial P$ in $P$).  

\vskip .1in
 \noindent For such a
representative of
${\cal F}$, we finally consider its restriction to each annular neighborhood $A_i$ of $\partial _i$.  Choose another arc $a$  whose endpoints are
{\sl not} in the windows (and again such an arc is essentially uniquely determined up to isotopy rel windows from a framing of the normal bundle to
$\partial P$ in $P$ in the natural way); orient $a$ and each component arc of $(\cup \alpha )
\cap A_i$ from
$\partial P_1$ to $\partial P$, and let $t_i$ be the signed (weighted) intersection number of $a$ with the (weighted)
arc family $(\cup \alpha)\cap A_i$, for $i=1,2,3$.  

\vskip .1in

\centerline{\epsffile{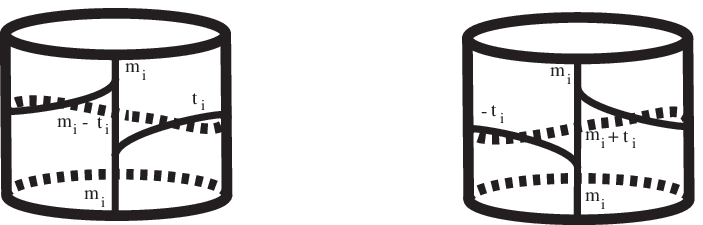}}

\vskip .1in

\hskip .6in{\bf 4a}~Right twisting for $t_i\geq 0$\hskip .3in{\bf 4b}~Left twisting for $t_i\leq 0$

\vskip .1in

\centerline{{\bf Figure 4}~Windowed measured foliations in the annulus}

\vskip .1in

\noindent As illustrated in Figure~4, all possible real twisting numbers $-m_i\leq t_i\leq
m_i$ arise provided
$m_i\neq 0$,  where again in this figure, the indicated ``weight'' of a component arc is the width of a band of
leaves parallel to the arc.  By performing Dehn twists along the core of the annulus, it likewise follows that every real twisting number
$t_i$ occurs provided $m_i\neq 0$.  Again, elementary topological considerations show that each windowed isotopy class of a windowed
measured foliation is uniquely determined by its invariants:

\vskip .2in

\noindent {\bf Lemma 1}~\it Points of ${\cal MF}(P)$
are uniquely determined by the triple $(m_i,t_i)\in{\bf R}_{\geq 0}\times{\bf R}$, which are subject
only to the constraint that $\forall i=1,2,3 (m_i=0\Rightarrow t_i\geq 0)$.  \rm

\vskip .1in

\vskip .2in

\noindent One difference between the Dehn-Thurston Lemma and Lemma~1 is that closed leaves are permitted in the
latter (but not in the former), where the coordinates $m_i=0$ and $t_i=|t_i|> 0$ correspond to the class of a
foliated annulus of width $t_i$ whose leaves are parallel to $\partial _i$. In the topology of projective
measured foliations, extensive twisting to the right or left about $\partial _i$ approaches the curve parallel
to
$\partial _i$.  One imagines identifying in the natural way the ray $\{0\}\times{\bf R}_{\geq 0}$ with the ray  
$\{0\}\times{\bf R}_{\leq 0}$ in the half plane ${\bf R}_{\geq 0}\times{\bf R}$ and thinks therefore of $(m_i,t_i)$ as lying in the
following quotient homeomorph of ${\bf R}^2$: 
$$R~=~({\bf R}_{\geq 0}\times{\bf R})/{\rm antipodal~map}.$$
We shall also require the subspace 
$$Z~=~({\bf Z}_{\geq 0}\times{\bf Z})/{\rm antipodal~map},$$
which corresponds to the collection of all disjointly embedded weighted curves and arcs in $F$ with endpoints
in the windows.

\vskip .1in

\noindent Arguing as above with an annular neighborhood of a pants decomposition, one concludes:

\vskip .2in

\noindent {\bf Theorem 2}~[23,32]~~\it Given an isotopy class of pants decomposition ${\cal P}$ of $F=F_{g,r}^s$, where each pants curve is
framed, there is a homeomorphism between ${\cal MF}(F)$ and the space of all pairs $(m_i,t_i)\in R$ as $i$ ranges over the elements of
${\cal P}$. Likewise, there is a
homeomorphism between ${\cal MF}_0(F)$ and the space of all pairs $(m_i,t_i)\in R$ as $i$ ranges over the elements of ${\cal
P}-\partial F$.  In particular, ${\cal PF}_0(F)\approx S^{3g-3+r+s}$ and ${\cal PF}(F)\approx S^{3g-3+2r+s}$.\rm

\vskip .2in

\noindent There is the following ``standard problem'' about which not much is known (on the torus, it devolves to greatest common divisors,
cf. [6], and see [11] for genus two):

\vskip .2in

\noindent {\bf Problem 2}~Given a tuple $\times _{i=1}^N (m_i,t_i)\in Z^N$, give a tractable expression in terms of Dehn-Thurston
or other coordinates for the number of components of the corresponding weighted family of curves and arcs.

\vskip .2in

\noindent There is an algorithm which leads to a multiply weighted curve from an integral measure on a general train track akin to that
on the torus gotten by serially ``splitting'' the track, cf. [6], but we ask in Problem~2 for a more closed-form expression.  See also Problem~3.
A related problem which also seems challenging is to describe ${\cal A'}(F)$ 
or ${Arc}'(F)$ in Dehn-Thurston coordinates on ${\cal MF}(F)$.  

\vskip .1in
 \noindent This class of curve and arc component counting
problems might be
approachable using the quantum path ordering techniques of [5,6] or with standard fermionic
statistical physics [57,29].

\vskip .3in

\noindent {\bf 3. Mapping class action on Dehn-Thurston coordinates}

\vskip .2in

\noindent As already observed by Max Dehn [8] in the notation of Theorem~2, a Dehn twist on the $i^{\rm th}$
pants curve in a pants decomposition of $F=F_{g,r}^s$ acts linearly, leaving invariant all coordinates
$(m_j,t_j)$, for
$j
\neq i$, and sending 
$$\eqalign{ 
m_i&\mapsto m_i,\cr 
t_i&\mapsto t_i\pm m_i.\cr
}\leqno{(\dag )}$$
\vskip .1in

\noindent As proved by Hatcher-Thurston [13] using Cerf theory, the two {\it elementary transformations~{\rm or}~moves} illustrated in Figure 5
act transitively on the set of all pants decompositions of any surface $F$.  Thus, any
Dehn twist acts on coordinates by the conjugate of a linear map, where the conjugating transformation
is described by compositions of ``elementary transformations'' on Dehn-Thurston coordinates
corresponding to the elementary moves.  
More explicitly, it is not difficult [23] to choose a finite
collection of pants decompositions of $F$ whose union contains all the curves in
Lickorish's generating set [20] and calculate the several compositions of elementary moves
relating them.

\vskip .1in

\noindent In any case, the calculation of the action of Dehn twist generators for $MC(F)$ thus
devolves to that of the two elementary moves on Dehn-Thurston coordinates.  This
problem was suggested in [8], formalized in [38], and solved in [23] as follows.

\vskip .1in

\noindent Let $\vee$ and $\wedge$ respectively denote the
binary infimum and supremum.

\vskip .1in

\noindent Given an arc family in the pair of pants, we introduce the notation $\ell _{ij}$ for
the arc connecting boundary components $i$ and $j$, for $i,j=1,2,3$, where the windowed isotopy class 
of $\ell _{ii}$ depends
upon the choice of twisting conventions.  Given a weighted arc family in the pair of pants, the
respective weights
$\lambda _{ij}$, for $\{ i,j,k\}=\{ 1,2,3\}$, on these component arcs are given in terms of the intersection
numbers $m_1,m_2,m_3$ by the following formulas:
$$\eqalign{
2\lambda _{ii}&=(m_i-m_j-m_k)\vee 0\cr
2\lambda_{ij}&=(m_i+m_j-m_k)\vee 0,~{\rm for}~i\neq j,\cr
}$$
and the intersection numbers are in turn given by
$m_i=2\lambda _{ii}+\lambda_{ij}+\lambda_{ik}$.

\vskip .2in

\epsffile{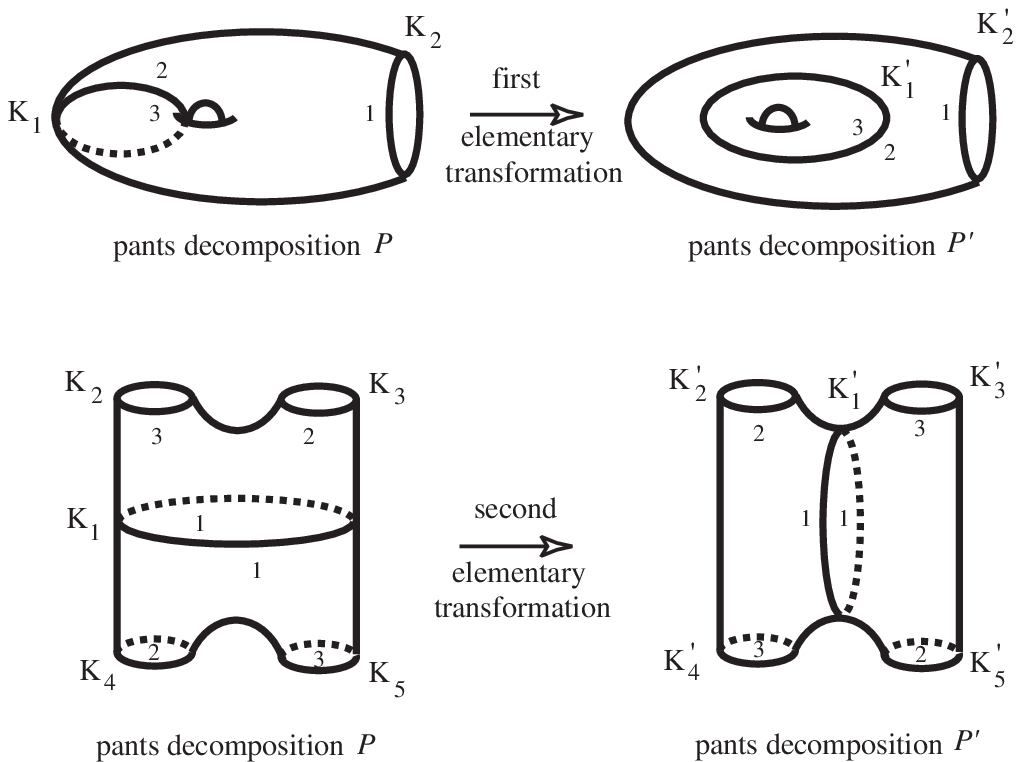}

\vskip .2in

\centerline{{\bf Figure 5}~The elementary transformations.}

\vskip .2in

\noindent {\bf Theorem 3}~~{\rm [23,32]}\it Adopt the standard twisting conventions, the enumeration of curves indicated
in Figure~5, and let $(m_i,t_i)$ denote the Dehn-Thurston coordinates of a measured foliation with
respect to the pants decomposition ${\cal P}$. 

\vskip .3in

\noindent{\bf First Elementary Transformation}~Let $\lambda _{ij}$ denote the weight of $\ell _{ij}$ with
respect to
${\cal P}$ and $\lambda '_{ij}$ the weight with respect to ${\cal P}'$, so in particular,
$r=\lambda_{12}=\lambda_{13}$.  Then the first
elementary transformation is given by the following formulas.
$$\eqalign{
\lambda _{11}'&=(r-|t_1|)\vee 0,\cr
\lambda' _{12}&=\lambda'_{13}=L+\lambda_{11},\cr
\lambda'_{23}&=|t_1|-L,\cr
t'_2&=t_2+\lambda _{11}+(L\wedge t_1)\vee 0,\cr
t_1'&=-sgn(t_1)~(\lambda_{23}+L),\cr
}$$
where $L=r-\lambda'_{11}$ and $sgn(x)\in\{\pm 1\}$ is the sign of $x\in{\bf R}$, with $sgn(0)=-1$.

\vskip .2in

\noindent {\rm (The formulas [24] above correct a typographical error in [32].)}

\vskip .2in

\noindent {\bf Second Elementary Transformation}~Let $\lambda_{ij}$ denote the weight of $\ell _{ij}$ in the
bottom pair of pants for ${\cal P}$ and $\kappa_{ij}$ the weight in the top pair or pants, and let
$\lambda _{ij}'$ denote the weight in the left pair of pants for ${\cal P}'$ and $\kappa '_{ij}$
in the right pair or pants.  The second elementary transformation is given by the following formulas.
$$\eqalign{
\kappa'_{11}&=\kappa_{22}+\lambda_{33}+(L-\kappa_{13})\vee 0+(-L-\lambda_{12})\vee 0,\cr
\kappa'_{22}&=(L\wedge\lambda_{11}\wedge(\kappa_{13}-\lambda_{12}-L))\vee 0,\cr
\kappa'_{33}&=(-L\wedge\kappa_{11}\wedge(\lambda_{12}-\kappa_{13}+L))\vee 0,\cr
\kappa'_{23}&=(\kappa_{13}\wedge\lambda_{12}\wedge(\kappa_{13}-L)\wedge(\lambda_{12}+L))\vee 0,\cr
\kappa'_{12}&=-2\kappa'_{22}-\kappa'_{23}+\kappa_{13}+\kappa_{23}+2\kappa_{33},\cr
\kappa'_{13}&=-2\kappa'_{33}-\kappa'_{23}+\lambda_{12}+\lambda_{23}+2\lambda_{22},\cr
\lambda'_{11}&=\lambda_{22}+\kappa_{33}+(K-\lambda_{13})\vee 0+(-K-\kappa_{12})\vee 0,\cr
\lambda'_{22}&=(K\wedge\kappa_{11}\wedge(\lambda_{13}-\kappa_{12}-K))\vee 0.\cr
\lambda'_{33}&=(-K\wedge\lambda_{11}\wedge(\kappa_{12}-\lambda_{13}+K))\vee 0,\cr
\lambda'_{23}&=(\lambda_{13}\wedge\kappa_{12}\wedge(\lambda_{13}-K)\wedge(\kappa_{12}+K))\vee 0,\cr
\lambda_{12}'&=-2\lambda_{22}'-\lambda'_{23}+\lambda_{13}+\lambda_{23}+2\lambda_{33},\cr
\lambda'_{13}&=-2\lambda'_{33}-\lambda'_{23}+\kappa_{12}+\kappa_{23}+2\kappa_{22},\cr
t_2'&=t_2+\lambda_{33}+((\lambda_{13}-\lambda'_{23}-2\lambda'_{22})\wedge(K+\lambda'_{33}-\lambda'_{22}))\vee
0,\cr
t_3'&=t_3-\kappa'_{33}+((L+\kappa'_{33}-\kappa'_{22})\vee
(\kappa'_{23}+2\kappa'_{33}-\lambda_{12}))\wedge 0,\cr
t_4'&=t_4-\lambda'_{33}+((K+\lambda'_{33}-\lambda'_{22})\vee (\lambda'_{23}+2\lambda'_{33}-\kappa_{12}))\wedge
0,\cr 
t_5'&=t_5+\kappa_{33}+((\kappa_{13}-\kappa'_{23}-2\kappa'_{22})\wedge(L+\kappa'_{33}-\kappa'_{22}))\vee
0,\cr
t_1'&=\kappa_{22}+\lambda_{22}+\kappa_{33}+\lambda_{33}-(\lambda'_{11}+\kappa'_{11}+(t_2'-t_2)+(t_5'-t_5))\cr
&~~~+~[(sgn(L+K+\lambda'_{33}-\lambda'_{22}+\kappa'_{33}-\kappa'_{22})]~(t_1+\lambda'_{33}+\kappa'_{33}),\cr
}$$ where $L=\lambda_{11}+t_1$, $K=\kappa_{11}+t_1$, and $sgn(x)\in\{\pm 1\}$ is the sign of $x\in{\bf R}$,
with 
$$sgn(0)=\cases{
+1,&if $\lambda_{12}-2\kappa'_{33}-\kappa'_{23}\neq 0$;\cr
-1,& otherwise.\cr}$$\rm

\vskip .2in

\noindent These formulas are derived in [23], in effect, by performing explicit isotopies of arcs in
certain covers of $F_{1,1}^0$ and $F_{0,4}^0$.  Their computer implementation has been useful
to some for analyzing specific mapping classes, e.g., fibered knot monodromies given by Dehn twists.

\vskip .1in

\noindent 
It is
notable that the formulas are piecewise-{\sl integral}-linear or PIL, cf. [37].  Furthermore, all of the
``corners in the PIL structure actually occur'', so the formulas are non-redundant in this sense; on the other
hand, a given word in Lickorish's generators (i.e., the composition of conjugates of linear mappings ($\dag$) by specific PIL transformations
given by finite compositions of the two elementary transformations) may not have ``all possible corners occur''. Insofar as
continuous concave PIL functions are in one-to-one correspondence with tropical polynomials [40,41], we are led to the following problem:

\vskip .2in

\noindent {\bf Problem 3}~Give a useful (piecewise) tropical description of the two elementary transformations.
One thus immediately derives a (piecewise) tropical polynomial representation of the mapping class groups.
What properties does it have, for instance under iteration?

\vskip .2in

\noindent As alternative coordinates, [9] describes a family of curves whose intersection numbers alone
coordinatize measured foliations of compact support (but there are relations), and presumably these intersection
numbers could be computed using Theorem 3.

\vskip .1in

\noindent We also wonder what are further applications or consequences of all these formulas.

\vskip .3in

\noindent{\bf 4. Pseudo-Anosov maps and the length spectrum of moduli space}

\vskip .2in

\noindent Thurston's original construction of pseudo-Anosov (pA) mappings [37] (cf. [9]) was generalized in
[23] (cf. [25,39]) to give the following recipe for their construction:

\vskip .2in

\noindent{\bf Theorem~4}~[23,25]~~\it Suppose that ${\cal C}$ and ${\cal D}$ are each families of
disjointly embedded essential simple closed curves so that each component of
$F-({\cal C}\cup{\cal D})$ is either disk, a once-punctured disk, or a boundary-parallel annulus.  Let $w$ be any word
consisting of Dehn twists to the right along elements of ${\cal C}$ and to the left along elements of ${\cal
D}$ so that the Dehn twist on each curve of ${\cal C}$ or ${\cal D}$ occurs at least once in $w$.
Then $w$ represents a pseudo-Anosov mapping class.\rm

\vskip .2in

\noindent {\bf Problem 4}~Does the recipe in Theorem~4 give virtually all pA maps?
That is, given a pA map $f$, is there some iterate $f^n$, for $n\geq 1$, so that $f^n$ arises from the recipe?

\vskip .1in

\noindent (This question from [23,25] is related [10] to the Ehrenpreis Conjecture, our Problem 14.) 

\vskip .1in

\noindent In relation to Problem~4, let us mention that there are still other descriptions of pA maps up to iteration, for instance
by Mosher [42] and in joint work of the author with Papadopoulos [45]; these descriptions are combinatorial rather than 
in terms of Dehn twists.

\vskip .1in

\noindent For a fixed surface $F$, consider the set of logarithms of dilatations of all pA maps
supported on
$F$.  This characteristic ``spectrum'' $\Sigma (F)\subseteq {\bf R}_{>0}$ of $F$ is precisely the Teichm\"uller
geoedsic length spectrum of Riemann's moduli space ${\cal M}(F)$.  The spectrum $\Sigma (F)$ is 
discrete.  (In fact, dilatations occur as spectral radii of integral-linear Perron-Frobenius
symplectomorphisms in a range of dimensions bounded above and below in terms of the topological type of
$F$.)

\vskip .1in

\noindent {\bf Problem 5}~For a given surface $F$, calculate $\Sigma (F)$.  More modestly, calculate the least
element of $\Sigma (F)$ or the least gap among elements of $\Sigma (F)$.  Characterize the number fields arising
as dilatations of pA maps on $F$.

\vskip .2in

\noindent Problems 4 and 5 are clearly related.  For instance, the recipe in Theorem~4 allows one to give
estimates on least elements in Problem~5, cf. [26,2].  McMullen [50] has also given estimates and examples.

\vskip .3in

\noindent {\bf 5. Arc complexes}

\vskip .2in

\noindent Refining the discussion in $\S$1, let $F_{g,\vec\delta}^s$ denote a bordered surface of type
$F_{g,r}^s$, with $r>0$ and $\vec\delta=(\delta _1,\ldots ,\delta _r )$ an $r$-dimensional vector of
natural numbers $\delta _i\geq 1$, where there are $\delta _i$ distinguished points on the
$i^{\rm th}$ boundary component of $F_{g,\vec\delta}^s$, for $i=1,\ldots, r$.  Construct an arc complex $Arc(F)$ as before as the
$PMC(F)$-orbits of isotopy classes of families of disjointly embedded essential and non-boundary parallel arcs connecting distinguished
points on the boundary.

\vskip .1in

\noindent Given two bordered surfaces $S_1,S_2$, we consider inclusions $S_1\subseteq S_2$, where the
distinguished points and punctures of $S_1$ map to those of $S_2$, and $S_1$ is a complementary component to an
arc family in $S_2$ (possibly an empty arc family if $S_1=S_2$).

\vskip .1in

\noindent Define the {\it type 1 surfaces} to be the following: $F_{1,(1,1)}^0$, $F_{0,(1,1,1,1)}^0$,
$F_{0,(1,1,1)}^1$, $F_{0,(1,1)}^2$.

\vskip .2in

\noindent {\bf Theorem 5}~[46]~\it ~~~The arc complex $Arc(F)$ of a
bordered surface $F$ is PL-homeomorphic to the sphere of dimension
$6g-7+3r+2s+\delta _1+\delta _2+\cdots +\delta _r$ if and only if 
$M\not\subseteq F$ for any type 1 surface $M$.  
In other words, $Arc(F)$ is spherical only
in the following cases:  polygons ($g=s=0$, $r=1$), multiply punctured polygons
($g=0$, $r=1$),
``generalized'' pairs of pants ($g=0$, $r+s=3$), the torus-minus-a-disk ($g=r=1$, $s=0$), and the once-punctured torus-minus-a-disk
($g=r=s=1$).  Only the type 1 surfaces have an arc complex which is 
a PL-manifold other than a sphere. \rm

\vskip .2in

\noindent {\bf Problem 6}~Calculate the topological type (PL-homeomorphism, homotopy, homology... type)
of the $Arc$-complexes.

\vskip .2in

\noindent The first non-trivial case is the calculation of the topological type of the PL-manifolds $Arc(M)$ for the four type 1
surfaces $M$.

\vskip .1in

\noindent Arc complexes as stratified spaces conjecturally have specific
singularities and topology described recursively as follows.  
A PL sphere is a type zero space.  A closed, connected, and simply connected manifold is
a type one space provided it occurs among a list of four specific such (non-spherical) manifolds of respective
dimensions 5,7,7, and 9, namely, the arc complexes of the four type one surfaces.  For $n>1$, define a type $n$ space to be
a finite polyhedron, defined up to PL-isomorphism, so that the link of each vertex in any compatible
triangulation is PL isomorphic to an iterated suspension of the join of at most two spaces of type less than
$n$.  

\vskip .1in

\noindent By Theorem 5, many links of simplices are indeed of this type, and we conjecture that any arc complex
is of some finite type.  (As explained in [46], a specific collapsing argument in the ``calculus of mapping cylinders'' 
would give a proof a of this conjecture.)

\vskip .1in

\noindent 
The non-Hausdorff space ${\cal PF}(F)/PMC(F)$ thus contains the stratified space
$Arc(F)$ as an open dense subset, explaining one classical (i.e., non-quantum)
aspect to Problem 1.  In light of the stratification of
$Arc$-complexes in general, one might hope to apply techniques such as [1,34] to address
parts of Problem 6.

\vskip .2in

\noindent {\bf Problem 7}~Devise a matrix model (cf. [28]) for the calculation of the Euler characteristics of $Arc$-complexes.

\vskip .1in

\noindent {\bf Problem 8}~[Contributed by the referee]~Does Theorem~5 say anything about the structure of the end of Riemann's moduli
space?  For instance, what is the homology near the end?

\vskip .2in

\noindent Take the one-point compactification $F^\times$ of $F=F_{g,r}^s$, where all of the $s\geq 0$ punctures of $F$ are identified to a single
point in $F^\times$, so $F^\times=F$ if and only if $s\leq 1$. 

\vskip .1in

\noindent Let
$\partial
$ denote the boundary mapping of the chain complex
$\{ C_p (Arc):p\geq 0\}
$ of  $Arc=Arc(F)$.   Suppose that
$\alpha$ is an arc family in
$F$ with corresponding cell $\sigma [\alpha ]\in C_p(Arc)$.  A codimension-one face of
$\sigma [\alpha ]$ of course corresponds to removing one arc from $\alpha$, and
there is a dichotomy on such faces $\sigma [\beta ]$ depending upon whether the rank
of the first homology of 
$F^\times-\cup\beta$ agrees with  or differs by one from that of $F^\times-\cup\alpha$.
This dichotomy decomposes $\partial $ into the sum of two operators $\partial =\partial
^1+\partial ^2$, where $\partial ^2$ corresponds to the latter case.

\vskip .1in

\noindent The operators
$\partial ^1,\partial ^2$ are a pair of anti-commuting differentials,
so there is a  spectral sequence
converging to $H_*(Arc)$ corresponding to the bi-grading
$$E^0_{u,v}=\{ {\rm chains~on}~\sigma [\alpha ]\in C_p(Arc):v=-{\rm rank}(H_1(F_\alpha))~{\rm
and}~u=p-v\} ,$$
where $\partial _1:E_{u,v}^0\to E_{u-1,v}^0$ and the differential of the $E^0$ term is $\partial _2:E_{u,v}^0\to
E_{u,v-1}^0$. 

\vskip .1in

\noindent It is not quite fair to call it a problem, nor a theorem since the argument is complicated and has not been independently checked, but
we believe that this spectral sequence collapses in its $E^1$-term to its top horizontal row except in dimension zero.  
Thus, the homology of $Arc$
is the $\partial _1$-homology of the $\partial _2$-kernels in the top row, and on the other hand, it follows from [31] that the 
$\partial _1$-homology of the top row itself agrees with that of uncompactified Riemann's moduli space ${\cal M}(F)$

\vskip .1in

\noindent As discussed in [30,46], the stratified structure of the arc complexes for bordered surfaces gives a corresponding
stratified structure to ${\cal A}(F)$ for punctured $F$.  This may be enough to re-visit the calculations of [19] and
[16,22] with an eye towards avoiding technical difficulties with the Deligne-Mumford 
compactification $\bar{\cal
M}(F)$.  

\vskip .2in

\noindent {\bf 6. Cell decompositions of ${\cal M}(F)$ and $\bar {\cal M}(F)$}

\vskip .2in

\noindent For the next several sections unless otherwise explicitly stated, surfaces $F$ will be taken to
be once-punctured and without boundary.  This is done for simplicity in order for instance that the moduli space ${\cal M}(F)$ itself, rather
than some ``decoration'' of it,  comes equipped with an ideal cell decomposition.  Nevertheless, the
theory extends, the discussion applies, and the problems and conjectures we articulate are intended in the
more general setting (of multiply punctured surfaces with a distinguished puncture and no boundary).

\vskip .1in

\noindent The basic combinatorial tool for studying moduli space ${\cal M}(F)$ is the $MC(F)$-invariant ideal cell
decomposition of the Teichm\"uller space ${\cal T}(F)$, and there are two effective constructions: from
combinatorics to conformal geometry using Strebel coordinates on fatgraphs [36,12,15], and from hyperbolic geometry to
combinatorics using the convex hull construction and simplicial coordinates [27].   (See [31] for further details.)  

\vskip .2in

\noindent{\bf Problem 9} [Bounded Distortion Conjecture] ~Given a hyperbolic structure on $F$,
associate its combinatorial invariant, namely, an ideal cell decomposition of $F$ together with the projective simplicial
coordinate assigned to each edge.  Take these projective simplicial coordinates as Strebel coordinates on the dual
fatgraph to build a conformal structure on
$F$.  The underlying map on Teichm\"uller space is of bounded distortion in the Teichm\"uller metric.

\vskip .2in

\noindent As posed by Ed Witten to the author in the early
1990's, a compelling problem at that time was to find an orbifold compactification of ${\cal M}(F)$ which
comes equipped with a cellular description in terms of suitably generalized fatgraphs.  Calculations such 
as [16,19,22] and more might then be performed using matrix
models derived from the combinatorics of this putative compactification.  
Perhaps the desired
compactification was the Deligne-Mumford compactification or perhaps another one.  The combinatorial
compactification of the previous section fails to provide an orbifold but rather another stratified
generalization of manifold.  

\vskip .1in

\noindent Guidance from Dennis Sullivan has recently led to the following
solution:

\vskip .2in

\noindent {\bf Theorem 6} \it Suppose that $F$ has only one puncture.  Then
$\bar {\cal M}(F)$ is homeomorphic to the geometric realization of the partially ordered set of $MC(F)$-orbits
of pairs 
$(\alpha, {\cal A})$, where $\alpha$ fills $F$, and the {\rm screen} ${\cal A}$
is a collection of subsets of $\alpha$ so that:

\vskip .1in

\leftskip .3in

\noindent {\rm [Fulton-MacPherson nest condition]} for any two $A,B\in {\cal A}$ which are not disjoint,
either
$A\subseteq B$ or $B\subseteq A$; 

\vskip .1in

\noindent {\rm [Properness]} $\cup{\cal A}$ is a proper subset of $\alpha$, and furthermore for any $A\in{\cal
A}$, we likewise have $\cup\{ B\in{\cal A}:A\neq B\subseteq A\}$ is a proper subset of $A$; 

\vskip .1in

\noindent {\rm [Recurrence]} for any
$a\in A\in{\cal A}$, there is an essential simple closed curve in $F$, meeting $\cup\alpha$ a minimal number
of times in its isotopy class, meeting only the arcs in $A$, and crossing $a$,  

\vskip .1in

\leftskip=0ex

\noindent where inclusion of ideal cell
decompositions induces the partial ordering on the set of pairs $(\alpha,{\cal A})$.\rm

\vskip .2in

\noindent In effect, the elements of the screen ${\cal A}$ detect how quickly hyperbolic lengths of arcs in
$\alpha$ diverge.  It is
a concise new description of a combinatorial structure  on $\bar {\cal M}(F)$, which bears similarities to
renormalization in physics, and it would be interesting to make this precise.  
Let us also remark that the proof of the previous theorem depends upon working in
the hyperbolic category, where the ``pinch curves'' in a nearly degenerate structure can be detected using
the coordinates.  (In fact, most of the required ideas and estimates are already described
in [30].)

\vskip .2in

\noindent {\bf Problem 10}~Though the (virtual) Euler characteristics are already known [28,56], devise a matrix model using
screens to calculate these invariants for ${\bar M}(F)$.

\vskip .2in

\noindent{\bf 7. Torelli groups}

\vskip .2in

\noindent Recall [52],[48] the Torelli group ${\cal I}_k(F)$, defined as
those mapping classes on $F$ fixing a basepoint (taken to be a puncture) that act trivially on the kth nilpotent quotient of the fundamental 
group of $F$.  Since the ideal cell decomposition of ${\cal T}(F)$ is invariant under $MC(F)$, it is in particular
invariant under each ${\cal I}_k(F)$, and so the quotient ``Torelli space'' $T_k(F)={\cal T}(F)/{\cal I}_k(F)$
is a manifold likewise admitting an ideal cell decomposition.

\vskip .1in

\noindent Recent work [48] with Shigeyuki Morita studies the
``Torelli tower'' 
$${\cal T}(F)\to\cdots\to T_{k+1}(F)\to T_k(F)\to\cdots\to T_1(F)\to {\cal M}(F)$$
of covers of Torelli spaces, each $T_k(F)$ a manifold covering the orbifold ${\cal M}(F)$.
In particular,
one essentially immediately derives infinite presentations for all of the higher Torelli groups as well as a
finite presentation for instance of the ``level $N$ (classical) Torelli groups'', i.e., the subgroup of
$MC(F)$ which acts identically on homology with ${\bf Z}/N$ coefficients.

\vskip .2in

\noindent{\bf Problem 11}~[Level $N$ Torelli Franchetta Problem]~What is the second cohomology group
of the level $N$ Torelli group?

\vskip .2in

\noindent The cell decomposition of $\bar {\cal M}(F)$ described before is compatible with the ideal cell
decompositions of Torelli spaces, i.e., the fatgraph dual to an ideal cell decomposition of $F$ admits a ``homology marking'' in the sense of
[48] as well as admitting the structure of screens.  There are thus  ``DM type'' boundaries of each Torelli space replete with an ideal cell
decomposition. It is natural to try to understand the topology of these DM-type bordifications of Torelli spaces to approach the following
class of problems:

\vskip .2in

\noindent {\bf Problem 12}~
Calculate various group-theoretic boundaries of mapping class and Torelli groups, for instance, Tits boundaries.

\vskip .2in

\noindent Further results in [48] arise by writing an explicit one cocyle representing the first Johnson
homomorphism.  In effect, one checks that the putative cycle represents a crossed homomorphism by verifying
the several essentially combinatorial constraints imposed by the cell decomposition, and then compares
with known values of the first Johnson homomorphism.  As discussed more fully in [48], one might
realistically hope to find similar canonical cocycles for the higher Johnson
homomorphisms as well.   In particular by work of Morita [53], the Casson invariant of a homology 3-sphere is
algorithmically calculable from the second Johnson homomorphism, so this is of special interest.

\vskip .1in

\noindent More explicitly, there is a combinatorial move, a ``Whitehead move'', that acts transitively on ideal triangulations of a
fixed surface (where one removes an edge $e$ from the triangulation to create a complementary quadrilateral
and then replaces $e$ with the other diagonal of this quadrilateral).  One seeks invariants of
sequences of Whitehead moves lying in an $MC(F)$-module satisfying three explicit combinatorial conditions
just as for the Johnson homomorphism.
There is thus a
kind of ``machine'' here for producing cocycles with values in various modules by solving for expressions
that satisfy certain explicit combinatorial constraints; several such invariant one cocycles have been
produced in this way on the computer but so far without success for constructing the
higher Johnson homomorphisms.

\vskip .1in

\noindent
Very recent work with Shigeyuki Morita and Alex Bene solves a related problem: the Magnus representations
(which are closely related to the Johnson homomorphisms [52]) in fact lift directly to the groupoid level of [48] as the Fox
Jacobians of appropriately enhanced Whitehead moves.  One works in the free fundamental group of the
punctured surface with basepoint distinct from the puncture, as discussed in the closing remarks of [48],
and the corresponding dual fatgraph comes equipped with a canonical maximal tree defined by greedily
adding edges to the tree while traversing a small circle around the puncture starting from the basepoint.  Perhaps these explicit
calculations of the Magnus representations might be of utility for instance to address the following standard
question:

\vskip .2in

\noindent {\bf Problem 13}~What are the kernels of the Magnus representations?

\vskip .2in

\noindent  Finally in [48] by taking contractions of powers of our canonical one cocycle, new
combinatorially explicit cycles and cocycles on ${\cal M}(F)$ are constructed which on the other hand generate the
tautological algebra.  It is natural to wonder about the extension of these classes to the screen model for
$\bar {\cal M}(F)$ and to the combinatorial compactification, and to revisit [16,19,22] in this context.

\vskip .4in

\noindent{\bf 8. Closed/open string theory}

\vskip .2in

\noindent This section describes recent work [49] with Ralph Kaufmann, where it turns out that the material of Section~2
together with a further combinatorial elaboration
describes  a reasonable model for the phenomenology of interactions of open and
closed strings.

\vskip .1in

\noindent Returning now to multiply punctured surfaces $F=F_{g,\vec\delta}^s$ with at least one distinguished point
in each boundary component as in Section~5, let $D$ denote the set
of distinguished points in the boundary, and let $S$ denote the set of punctures.  Fix a set
${\cal B}$ of ``brane labels'', let ${\cal P}({\cal B})$ denote its power set, and define a {\it brane labeling} to be a function $\beta:D\cup
S\to {\cal P}({\cal B})$ subject to the unique constraint that if the value $\emptyset\in{\cal P}({\cal
B})$ is taken at a point on the boundary of $F$, then this point is the only
distinguished point in its boundary component. Complementary components to the distinguished points in the
boundary are called ``windows''.  Given a brane labeling $\beta$ on a windowed surface $F$, define the set
$D(\beta)=\{ d\in D:\beta (d)\neq\emptyset\}$, and consider proper isotopy classes rel $D$ of arc families
with endpoints in windows, where the arcs are required to be non-boundary
parallel in $F-D(\beta )$.  Say that such an arc family is {\it exhaustive} if each window has at least one incident arc,
and construct
the space $\widetilde{Arc}(F,\beta)$ of positive real weights on exhaustive arc families.
\vskip .1in

\noindent There are various geometric operations on the spaces $\{ \widetilde{Arc}(F,\beta )\}$ induced by
gluing together measured foliations along windows provided the total weights agree on the windows to be
glued (taking unions of brane labels when combining distinguished points).  It is shown in [49] that these
geometric operations descend to the level of suitable chains on the spaces $\widetilde{Arc}(F,\beta) $ and
finally to the level of the integral homology groups of
$\widetilde{Arc}(F,\beta )$.  These algebraic operations on homology satisfy the expected
``operadic'' equations of open/closed string theory, and new equations can be discovered as well.

\vskip .1in

\noindent In effect, [49] gives the string field theory of a single point, and natural questions are already discussed in detail in [49]
including  the ``passage to conformal field theory'' which seems to be provided by corresponding operations not for exhaustive arc families but
rather for quasi filling arc families; an obvious challenge is to perform meaningful calculations in CFT.  Calculate the homology groups of
exhaustive or quasi filling arc families in brane labeled windowed surfaces, and calculate the homology groups of their combinatorial or
DM-type compactifications (the former problem already articulated as our Problem~6 and the latter surely
also posed elsewhere in this volume).  Organize and understand the many algebraic relations of open/closed
strings on the level of homology.
Find the BRST operator.  Introduce excited strings, i.e., extend to string field theory of realistic targets.
Is the Torelli structure of any physical significance?

\vskip .2in

\noindent {\bf 9. Punctured solenoid}

\vskip .2in

\noindent A problem well-known in the school around Dennis Sullivan is:

\vskip .1in

\noindent {\bf Problem 15}~[Ehrenpreis Conjecture]~Given two closed Riemann surfaces, there are finite unbranched covers
with homeomorphic total spaces which are arbitrarily close in the Teichm\"uller metric.

\vskip .1in

\noindent This section describes joint work with Dragomir \v Sari\' c on related universal constructions in Teichm\"uller theory [47].

\vskip .1in

\noindent As a tool for understanding dynamics and geometry, Sullivan defined the hyperbolic solenoid [33] as the inverse
limit of the system of finite-sheeted unbranched pointed covers of any fixed closed oriented surface of negative Euler characteristic.  Following 
Ahlfors-Bers theory, he developed its Teichm\"uller theory and studied the
natural dynamics and geometry, in particular introducing two principal mapping class groups,
the continuous ``full mapping class group'' and the countable ``baseleaf preserving mapping class group''.
He furthermore showed that the Ehrenpreis Conjecture is equivalent to the conjecture that the latter group has dense orbits in the Teichm\"uller
space of the solenoid.  [55] contains many basic results and open problems about Sullivan's solenoid. 

\vskip .1in

\noindent Following this general paradigm, one defines the punctured hyperbolic solenoid $S$ as the inverse limit over all finite-index subgroups
$\Gamma$ of $PSL_2({\bf Z})$ of the system of covers ${\cal U}/\Gamma$ over the modular curve ${\cal U}/PSL_2({\bf Z})$, where ${\cal U}$
denotes the upper half-space.  In effect, branching is now permitted but only at the missing punctures covering the three orbifold points of the
modular curve.  

\vskip .1in

\noindent Let us build a particular model space homeomorphic to the punctured solenoid.  Take $\Gamma = PSL_2({\bf Z})$ acting by fractional
linear transformations on ${\cal U}$, and let $\hat\Gamma$ denote its profinite compleion.  Thus, $\gamma\in\Gamma$ acts naturally on ${\cal
U}\times\hat\Gamma$ by $\gamma\cdot (z,t) =(\gamma z,t\gamma ^{-1})$, and the quotient is homeomorphic to $S$.

\vskip .1in

\noindent Truly the entire decorated Teichm\"uller theory of a punctured surface of finite type [27] extends appropriately to $S$: there are
global coordinates, there is an ideal cell decomposition of the decorated Teichm\"uller space, and there is an explicit non-degenerate two-form,
the latter of which are invariant under the action of the baseleaf preserving subgroup $Mod_{BLP}(S)$ of the full mapping class group
$Mod(S)$; generators for
$Mod_{BLP}(S)$ are provided by appropriate equivariant Whitehead moves, and a complete set of relations has recently
been derived as well in further recent work with Dragomir \v Sari\' c and Sylvain Bonnot.  

\vskip .1in

\noindent There are a number of standard questions (again see [55]):
Is $Mod_{BLP}(S)$ finitely generated?  Do the ``mapping class like'' elements generate $Mod_{BLP}(S)$?  It seems to be a deep question
which tesselations of the disk, other than the obvious so-called ``TLC'' ones, arise from the convex hull construction applied to the
decorated solenoid in [47].  Is there any relationship between $Mod(S)$ and the absolute Galois group? 
(This cuts both ways since only partial information is known about either group; notice that the action of the latter on $\hat{\Gamma}$ is
explicit as a subgroup of the Grothendieck-Teichm\"uller group [51].) 

\vfill\eject

\noindent{\bf Bibliography}

\vskip .1in

\noindent [1]~~~Nils Andreas Baas, ``Bordism theories with singularities'', Proceedings of the Advanced Study Institute on
Algebraic Topology (Aarhus Univ., Aarhus, 1970), Vol. I, pp. 1-16. Various Publ. Ser., No. 13, Mat. Inst., Aarhus Univ.,
Aarhus, 1970.

\vskip .1in

\noindent [2]~~~Max Bauer,
``An upper bound for the least dilatation'', {\it Trans. Amer. Math. Soc.} {\bf 330} (1992), 361-370. 

\vskip .1in

\noindent [3]~~~Moira Chas and Dennis P. Sullivan, ``String topology'', 
to appear {\it Ann. Math.}, preprint math.GT/9911159,

\vskip .1in

\noindent [4]~~~---, 
``Closed string operators in topology leading to Lie bialgebras and 
higher string algebra'',  The legacy of Niels Henrik Abel,  
Springer, Berlin (2004), 771-784, Preprint math.GT/0212358.

\vskip .1in

\noindent [5]~~~L.~Chekhov and V.~Fock, {``A quantum Techm\"uller space''}, {\it Theor. Math.
Phys.}, {\bf 120} (1999) 1245-1259; {``Quantum mapping class group,
pentagon relation, and geodesics''}
{\it Proc. Steklov Math. Inst.}
{\bf 226} (1999) 149-163.

\vskip .1in

\noindent [6]~~~Leonid Chekhov and R. C. Penner,
``On the quantization of
Teichm\"uller and Thurston theories'', to appear: Handbook of Teichm\"uller theory,
European Math Society, 
ed. A. Papadopoulos, 
math.AG/0403247.

\vskip .1in

\noindent [7]~~~Ralph L.\ Cohen and John D.S.\ Jones
``A homotopy theoretic realization of string topology'', {\it Math. Ann.} {\bf 324} (2002), 773-798.

\vskip .1in

\noindent [8]~~~Max Dehn, Lecture notes from Breslau, 1922, Archives of the University of Texas at Austin.

\vskip .1in

\noindent [9]~~A. Fathi, F. Laudenbach, V. Poenaru, ``Travaux de Thurston sur les Surfaces'',
{\it Asterisques}~{\bf 66-67} (1979), Sem. Orsay, Soc. Math. de France.

\vskip .1in

\noindent [10]~~
Tim Gendron, 
``The Ehrenpreis conjecture and the moduli-rigidity gap'',
Complex manifolds and hyperbolic geometry (Guanajuato, 2001), 
{\it Contemp. Math.} {\bf 311}, 207-229.

\vskip .1in

\noindent [11]~~Andrew Haas and Perry Susskind,
``The connectivity of multicurves determined by integral weight train tracks'', {\it Trans. Amer. Math.
Soc.} {\bf 329} (1992), 637-652.

\vskip .1in

\noindent [12]~~ John L. Harer, 
``The virtual cohomological dimension of the mapping class group of an 
orientable surface'',  {\it Invent. Math.}  {\bf 84}  (1986),  157-176.

\vskip .1in

\noindent [13]~~A. Hatcher and W. Thurston, ``A presentation for the mapping class group of a closed orientble surface'', {\it
Topology}~{\bf 19} (1980), 221-137.

\vskip .1in
                                                                                
\noindent [14]~~William J. Harvey,
``Boundary structure of the modular group'',
Riemann surfaces and related topics: Proceedings of the 1978 Stony Brook
Conference,  {\it Ann. of Math. Stud.} {\bf 97},
Princeton Univ. Press, Princeton, N.J., (1981), 245-251.

\vskip .1in

\noindent [15]~~John H. Hubbard and Howard Masur,
``Quadratic differentials and foliations'',  {\it Acta Math.}  {\bf 142}
(1979), 221-274.

\vskip .1in

\noindent [16]~~Kyoshi Igusa,
``Combinatorial Miller-Morita-Mumford classes and Witten cycles'',
{\it  Algebr. Geom. Topol.}  {\bf 4}  (2004), 473-520.

\vskip .1in

\noindent [17]~~Rinat M. Kashaev,
{``Quantization of Teichm\"uller spaces and the
quantum dilogarithm''}, {\it Lett. Math. Phys.}, {\bf 43}, 
(1998), 105-115; q-alg/9705021.

\vskip .1in

\noindent [18]~~Ralph Kaufman, Muriel Livernet, R. C. Penner, ``Arc operads and arc
algebras'', {\it Geom. Topol.} {\bf 7} (2003), 511-568.

\vskip .1in

\noindent [19]~~Maxim Kontsevich,
``Intersection theory on the moduli space of curves and the
matrix Airy function'',  {\it Comm. Math. Phys.}  {\bf  147}  (1992), 1-23.

\vskip .1in

\vskip .1in

\noindent [20]~~W. B. R. Lickorish, ``A finite set of generators for the homeotopy group of a two-manifold'' {\it Proc. Camb.
Phil. Soc.} {\bf 60} (1964), 769-778; Corrigendum,
{\it Proc. Camb. Phil. Soc.} {\bf 62} (1966), 679-681.

\vskip .1in

\noindent [21]~~Yuri I. Manin, 
``Von Zahlen und Figuren'',
Talk at the International Conference ``GŽomŽtrie au vingtime cicle: 1930--2000'', Paris, Institut Henri PoincarŽ, Sept. 2001 
27 pp., preprint math.AG/0201005.

\vskip .1in

\noindent [22]~~Gabriel Mondello,
``Combinatorial classes
on $\overline{M}\sb {g,n}$ are tautological'',  {\it Int. Math. Res. Not.}
(2004), 2329-2390.

\vskip .1in

\noindent [23]~~R. C. Penner, ``The action of the mapping class group on isotopy classes of curves and arcs in 
surfaces",  thesis,  Massachusetts Institute of Technology (1982), 180 pages.

\vskip .1in

\noindent [24]~~---,
``The action of the mapping class group on curves in surfaces",  
{\it L'Enseignement Mathematique }{\bf 30} (1984),  39-55.

\vskip .1in

\noindent [25]~~---,
``A construction of pseudo-Anosov homeomorphisms", 
{\it Proceedings of the American Math Society} {\bf 104} (1988), 1-19.

\vskip .1in

\noindent [26]~~---,
``Bounds on least dilatations",  {\it Trans. Amer. Math. Soc.}   
{\bf 113} (1991), 443-450.

\vskip .1in

\noindent [27]~~---,
``The decorated Teichm\"uller space of  punctured surfaces", 
{\it Communications in Mathematical  Physics}  {\bf 113}   (1987),  299-339.

\vskip .1in

\noindent [28]~~---,
``Perturbative series and the moduli space of Riemann surfaces",   
{\it Journal  of  Differential Geometry} {\bf  27}  (1988),  35-53.

\vskip .1in

\noindent [29]~~~,
``An arithmetic problem in surface geometry'', {\it
The Moduli Space of Curves}, Birkh\"auser (1995), eds. R. Dijgraaf, 
C. Faber, G. van der Geer, 427-466.

\vskip .1in

\noindent [30]~~---,
``The simplicial compactification of Riemann's moduli space'', 
Proceedings of the 37th Taniguchi Symposium, World Scientific (1996), 237-252.

\vskip .1in

\noindent [31]~~---,
``Decorated Teichm\"uller space of bordered
surfaces'',  {\it Communications in Analysis and  Geometry} {\bf 12} (2004), 793-820, math.GT/0210326.

\vskip .1in

\noindent [32]~~
R. C. Penner with John L. Harer  {\it Combinatorics of Train Tracks},
Annals of Mathematical Studies  {\bf 125},  Princeton Univ. Press (1992);
second printing (2001).

\vskip .1in

\noindent [33]~~Dennis Sullivan, {\it Linking the universalities of
Milnor-Thurston, Feigenbaum and Ahlfors-Bers}, Milnor Festschrift,
Topological methods in modern mathematics (L. Goldberg and A.
Phillips, eds.), Publish or Perish, 1993, 543-563.
preprint.

\vskip .1in

\noindent [34]~~---,
``Combinatorial invariants of analytic spaces'', Proc. of Liverpool Singularities--Symposium, I (1969/70)
Springer, Berlin (1971), 165-168.  

\vskip .1in

\noindent [35]~~Joerg Teschner
``An analog of a modular functor from quantized TeichmŸller theory'', to appear
Handbook of Teichm\"uller theory, European Math Society (ed. A. Papadopoulos),
math.QA/0405332. 

\vskip .1in
                                                                                
\noindent [36]~~Kurt Strebel, Quadratic Differentials, {\it Ergebnisse
der Math. und ihrer Grenzgebiete}, Springer-Verlag, Berlin (1984).

\vskip .1in

\noindent [37]~~William P. Thurston, ``On the geometry and dynamics of diffeomorphisms of surfaces'',
{\sl Bull. Amer. Math. Soc.}, {\bf 19} (1988) 417--431.

\vskip .1in

\noindent [38]~~---, ``Three dimensional manifolds, Kleinian groups and hyperbolic geometry'', 
{\it Bull. Amer. Math. Soc.}~{\bf 6}(1986),
357-381.

\vskip .1in

\noindent [39]~~Albert Fathi,
``DŽmonstration d'un thŽorme de Penner sur la composition des twists de Dehn'', {\it Bull. Soc. Math. France} {\bf 120} (1992),
467-484.

\vskip .1in

\noindent [40]~~ Anatol N. Kirillov, 
``Introduction to tropical combinatorics'', 
Physics and combinatorics, 2000 (Nagoya), 82-150, 
World Sci. Publishing, River Edge, NJ, 2001. 

\vskip .1in

\noindent [41]~~J\"urg Richter-Gebert, Bernd Sturmfels, Thorsten Theobald, ``First steps in tropical geometry'', to appear in Proc. of
COnf. on Idempotent Mathemaicals and Mathematical Physics, {\it Contemp. Math.}, math.AG/0306366.

\vskip .1in

\noindent [42]~~Lee Mosher, ``The classification of pseudo-Anosovs"  in {\it Low-dimensional topology and Kleinian groups} (Coventry/Durham,
1984), 13-75, {\it London Math. Soc. Lecture Note Ser.}, {\bf 112}, Cambridge Univ. Press, Cambridge, 1986. 

\vskip .1in

\noindent [43]~~William A. Veech, ``Dynamics over TeichmŸller space'', {\it  Bull. Amer. Math. Soc. }~{\bf 14} (1986), 103-106.

\vskip .1in

\noindent [44]~~Howard Masur, ``Ergodic actions of the mapping class group'', {\it Proc. Amer. Math. Soc}~{\bf 94} (1985), 455-459.

\vskip .1in

\noindent [45]~~Athanase Papadopoulos and R. C. Penner, ``A construction of pseudo-Anosov homeomorphisms", 
{\it Proc. Amer. Math Soc.} {\bf 104} (1988), 1-19.

\vskip .1in

\noindent [46] R. C. Penner, ``The structure and singularities of arc complexes'', preprint (2004) math.GT/0410603.

\vskip .1in

\noindent [47]~~R. C. Penner and Dragomir \v Sari\' c, ``Teichm\"uller theory of the punctured solenoid'', 
preprint (2005), math.DS/0508476.

\vskip .1in

\noindent [48]~~S. Morita and R. C. Penner, ``Torelli groups, extended Johnson homomorphisms, and new cycles on the moduli space of curves'', 
preprint (2006), math.GT/0602461.

\vskip .1in

\noindent [49]~~Ralph M. Kaufmann and R. C. Penner, ``Closed/open string diagrammatics'' (2006), to appear in 
{\it Nucl. Phys B}, math.GT/0603485.

\vskip .1in

\noindent [50]
~~Curt McMullen, ``Billiards and Teichm\"uller curves on Hilbert modular surfaces'', {\it J. Amer. Math. Soc.} {\bf 16} (2003), 857Ð885. 

\vskip .1in

\noindent [51]~~Leila Schneps, ``The Grothendieck-Teichm\"uller group $\widehat{GT}$: a survey'', Geometric Galois Actions (L. Schneps and
P. Lochak, eds.) London Math Society Lecture Notes {\bf 242} (1997), 183-203.

\vskip .1in

\noindent [52]~~Shigeyuki Morita,
``{ Abelian quotients of subgroups of the mapping class
     group of surfaces}'',
{\it Duke Math. J.} {\bf 70}
(1993), 699--726.

\vskip .1in

\noindent [53]~~Shigeyuki Morita, ``Casson's invariant for homology 3-spheres and characteristic classes of surface bundles I'' 
{\it Topology} {\bf 28} (1989),
305-323.

\vskip .1in

\noindent [54]~~Yuri I. Manin,
``Real multiplication and noncommutative
geometry'', {\it The legacy of Niels Henrik Abel},
eds. O.~A.~Laudal and R.~Piene, Springer Verlag,
Berlin 2004, 685-727, preprint math.AG/0202109.

\vskip .1in

\noindent [55]~~Chris Odden, ``The baseleaf preserving mapping class group of the universl hyperbolic solenoid'',
{\it Trans. A.M.S.} {\bf 357} (2004), 1829-1858.

\vskip .1in

\noindent [56]~~John L. Harer and Don Zagier, 
``The Euler characteristic of the moduli space of curves'',
{\it Invent. Math.} {\bf 85} (1986), 457-485.

\vskip .1in

\noindent [57]~~Claude Itzykson and Jean-Michel Drouffe, Statistical Field Theory, volume 1,
Cambridge Monographs on Mathematical Physics, Cambridge Univ. Press (1989).

\vfill\eject

\bye